\documentclass{article}    

\usepackage{amsmath}
\usepackage{amstext}
\usepackage{amsthm}
\usepackage{amsfonts}

\usepackage[dvips]{graphicx}

\theoremstyle{plain}

\newtheorem{corollary}{Corollary}

\newtheorem{lemma}{Lemma}

\newtheorem{proposition}{Proposition}
\newtheorem{remark}{Remark}

\newtheorem{theorem}{Theorem}
\numberwithin{equation}{section}

\begin{document}

\title{ On the volume of spherical Lambert cube
\thanks{Partially supported by Russian Foundation for
Basic Research (Grant 00-15-96165)}}
\author{Dmitriy Derevnin
and Alexander Mednykh }

\maketitle

\begin{abstract}
{The calculation of volumes of polyhedra in the three-dimensional
Euclidean, spherical and hyperbolic spaces is very old and
difficult problem. In particular, an elementary formula for volume
of non-euclidean simplex is
still unknown. One of the simplest polyhedra is the Lambert cube $Q({\alpha }%
,\beta ,\gamma )$. By definition, $Q({\alpha },\beta ,\gamma )$ is
a combinatorial cube, with dihedral angles ${\alpha },\beta $ and
$\gamma $ assigned to the three mutually non-coplanar edges and
right angles to the remaining. The hyperbolic volume of Lambert
cube was found by Ruth Kellerhals (1989) in terms of the
Lobachevsky function $\Lambda (x)$. In the present paper the
spherical volume of $Q({\alpha },\beta ,\gamma )$ is defined in
the terms of the function
\begin{equation*}
\delta \left( \alpha ,\theta \right) =\int\limits_{\theta }^{\frac{\pi }{2}}{%
\log }\left( {1-\cos 2\alpha \cos 2\tau }\right) \frac{d\tau
}{\cos 2\tau }
\end{equation*}
which can be considered as a spherical analog of the function
\begin{equation*}
\Delta \left( \alpha ,\theta \right) =\Lambda \left( \alpha
+\theta \right) -\Lambda \left( \alpha -\theta \right) .
\end{equation*}
}


\end{abstract}

\section{{{Introduction}}}

The calculation of volumes of polyhedra in Euclidean, spherical
and hyperbolic spaces is very old and difficult problem. The main
principles for volume calculations in non-euclidean geometries
were given in 1836 by Lobachevsky \cite{L} and in 1852 by
Schl\"{a}fli \cite{Sh}. In particular, they have found the volumes
of the orthogonal three-dimensional simlicies (orthoschemes).

In general, every polyhedron can be decomposed into a finite number of
orthoschemes. But, in spite of this, an elementary formula for volume of
non-euclidean simplex is still unknown\footnote{%
Just recently, such a formula was obtained in \cite{CK} and
\cite{MY}. Simple proof of the formula can be found in \cite{U}.}.

About 1935 Coxeter \cite{C} revived interest in the work of these
two authors by developing an integration method for non-euclidean
orthoschemes of dimension three. This method was generalized by
B\"{o}hm \cite{B} for spaces of constant non-vanishing curvature
of arbitrary dimension. Further advance in this direction was
achieved in the papers by Vinberg \cite{V}, Milnor \cite{M} and
Ruth Kellerhals \cite{K}. The simplest generalization of notion of
three-dimensional orthoscheme is
the Lambert cube $Q({\alpha },\beta ,\gamma )$. Recall that $Q({\alpha }%
,\beta ,\gamma )$ is a combinatorial cube, with dihedral angles ${\alpha }%
,\beta $ and $\gamma $ assigned to the three mutually non-coplanar
edges and right angles to the remaining.

The Lambert cube $Q({\alpha},\beta,\gamma)$ can be realized in the
hyperbolic space (\cite{{K},{HLM}}) if $0<{\alpha},\beta,\gamma<\frac{\pi%
}{2}$ and in the spherical space \cite{D} if $\frac{\pi}{2}<{\alpha}%
,\beta,\gamma<\pi$.

The hyperbolic volume of $Q({\alpha},\beta,\gamma)$ was found in \cite{K} in
terms of the Lobachevsky function
\begin{equation*}
\Lambda(x)=-\int\limits_{0}^{x}{\log}\left| {2\sin t}\right| dt.
\end{equation*}

In the present paper the spherical volume of $Q({\alpha },\beta
,\gamma )$ is defined in the terms of the function
\begin{equation*}
\delta \left( \alpha ,\theta \right) =\int\limits_{\theta }^{\frac{\pi }{2}}{%
\log }\left( {1-\cos 2\alpha \cos 2\tau }\right) \frac{d\tau }{\cos 2\tau }
\end{equation*}
which can be considered as a spherical analog of the function
\begin{equation*}
\Delta \left( \alpha ,\theta \right) =\Lambda \left( \alpha +\theta \right)
-\Lambda \left( \alpha -\theta \right) \text{.}
\end{equation*}

The main result of the work is the following theorem.
\begin{theorem}
\label{mm}{{{The volume of a spherical Lambert cube $Q(\alpha ,\beta ,\gamma
)$ with essential angles} ${\alpha },\beta $ and $\gamma $, $\frac{\pi }{2}%
<\alpha ,\beta ,\gamma <\pi $ {is given by the formula }
\begin{equation*}
V(\alpha ,\beta ,\gamma )=\frac{1}{4}\left( \delta \left( \alpha ,\theta
\right) +\delta \left( \beta ,\theta \right) +\delta \left( \gamma ,\theta
\right) -2\delta \left( \frac{\pi }{2},\theta \right) -\delta \left( {0}%
,\theta \right) \right) ,
\end{equation*}
{where}
\begin{equation*}
\delta \left( \alpha ,\theta \right) =\int\limits_{\theta }^{\frac{\pi }{2}}{%
\log }\left( {1-\cos 2\alpha \cos 2\tau }\right) \frac{d\tau }{\cos 2\tau }
\end{equation*}
{and} $\theta $, $\frac{\pi }{2}<\theta <\pi $ {is the principal
parameter defined by}
\begin{align*}
\tan ^{2}\theta & =-p+\sqrt{p^{2}+L^{2}M^{2}N^{2}}, \\
p& =\frac{L^{2}+M^{2}+N^{2}+1}{2}, \\
L& =\tan \alpha ,M=\tan \beta ,N=\tan \gamma .
\end{align*}
}}
\end{theorem}

{It is interesting to compare properties of the function $\delta (\alpha
,\theta )$ with the properties of the function }$\Delta \left( \alpha
,\theta \right) $.{{\ }Notice that the hyperbolic volumes of knots,
orbifolds and cone manifolds in many cases can be expressed in terms of $%
\Delta \left( \alpha ,\theta \right) $ (see
\cite{{Th},{MV},{V},{K}}). In particular, the main result of the
work \cite{K} can be rewrited in the following form. }

\begin{theorem}
The volume of a hyperbolic Lambert cube $Q(\alpha ,\beta ,\gamma )$ with
essential angles ${\alpha },\beta $ and $\gamma $ , $0<\alpha ,\beta ,\gamma
<\frac{\pi }{2}$ is given by the formula
\begin{equation*}
V(\alpha ,\beta ,\gamma )=\frac{1}{4}\left( \Delta \left( \alpha ,\theta
\right) +\Delta \left( \beta ,\theta \right) +\Delta \left( \gamma ,\theta
\right) -2\Delta \left( \frac{\pi }{2},\theta \right) -\Delta \left(
0,\theta \right) \right) ,
\end{equation*}
where $\theta $, $0<\theta <\frac{\pi }{2}$ is the principal
parameter defined by conditions
\begin{align*}
\tan ^{2}\theta =p+\sqrt{p^{2}+L^{2}M^{2}N^{2}}, \\
p=\frac{L^{2}+M^{2}+N^{2}+1}{2}, \\
L=\tan \alpha ,M=\tan \beta ,N=\tan \gamma .
\end{align*}
\end{theorem}

\section{Metric properties of Lambert cube}

{Consider the Euclidean space }$\mathbb{R}${$^{4}=\left( R^{4},d\sigma
^{2}\right) $ with the metric }

{\
\begin{equation*}
d\sigma^{2}=\frac{{dx}^{2}}{A^{2}}+\frac{{dy}^{2}}{B^{2}}+\frac{{dz}^{2}}{%
C^{2}}+{dt}^{2}
\end{equation*}
induced by the scalar product
\begin{equation*}
((x,y,z,t),(x^{\prime},y^{\prime},z^{\prime},t^{\prime}))=\frac{xx^{\prime}}{%
A^{2}}+\frac{yy^{\prime}}{B^{2}}+\frac{zz^{\prime}}{C^{2}}+tt^{\prime},
\end{equation*}
where $A$,$B$ and $C$ are given positive numbers. }

{Define
\begin{equation*}
S^{3}=\{v=(x,y,z,t){\in}\mathbb{R}^{4}:(v,v)=1\}.
\end{equation*}
The sphere $S^{3}$ endowed with the metric $d\sigma^{2}$ becomes the
spherical space
\begin{equation*}
\mathbb{S}^{3}=({S}^{3},d\sigma^{2})
\end{equation*}
with the constant Gaussian curvature $k=1$. To see that we set }$ds^{2}={dx}%
^{2}+{dy}^{2}+{dz}^{2}+{dt}^{2}$.{\ Then the mapping $\left(
R^{4},ds^{2}\right) \longrightarrow$ $\left(
R^{4},d\sigma^{2}\right) $ defined by
}$(x,y,z,t)\rightarrow(Ax,By,Cz,t)$ is an isometry sending the
unit sphere in {$\left( R^{4},ds^{2}\right) $ onto
}$\mathbb{S}^{3}$.

{At the same time we consider }$R${$^{4}$ as a projective space }$%
RP ${$^{3}$ with the homogeneous coordinates $(x:y:z:t)$. }

{As the Klein model }$\mathbb{K}${\ of the spherical space }$\mathbb{S}${$%
^{3}$ choose the hyperplane }$\mathbb{K}${$=\{(x,y,z,t){\in }\mathbb{R}%
^{4}:t=1\}$, which will be identified further with the Euclidean space }$%
\mathbb{R}${$^{3}=\{(x,y,z):x,y,z{\in }R\}$. }

{There exists one-to-one correspondence between }$\mathbb{K}${\ and the
upper semisphere }$\mathbb{S}${$_{+}^{3}=\{(x,y,z,t){\in }\mathbb{S}%
^{3}:t>0\}$ formed by projection from the origin. }

{Reflections in coordinate planes and rotations in coordinate lines of }$%
\mathbb{R}${$^{4}$ are both Euclidean and spherical isometries. }

\begin{figure}
\resizebox{0.95\hsize}{!}{\includegraphics*{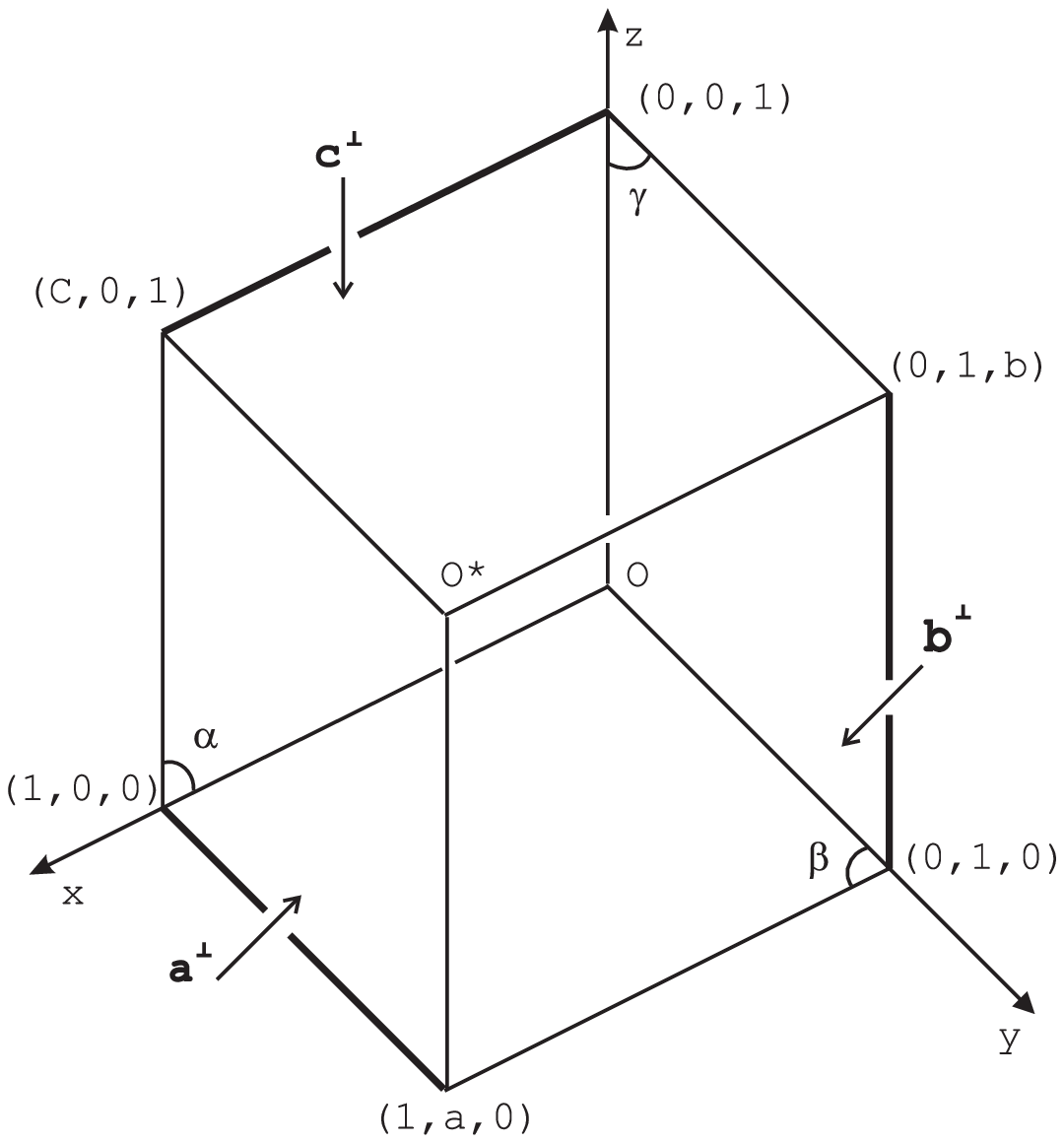}}
\caption{} \label{fig1}
\end{figure}

{Following to \cite{HLM} we realize the spherical Lambert cube $Q({\alpha }%
,\beta,\gamma)$ as a projection of the Euclidean polyhedron
$P(a,b,c)$ represented on the Fig. \ref{fig1} . }

{We can assume that the essential angles ${\alpha},\beta$ and
$\gamma$ are formed by the pairs of the planes
$\{a^{\perp},z=0\}$, $\{b^{\perp},x=0\} $ and $\{c^{\perp},y=0\}$.
The projective equations of the planes $a^{\perp}$, $b^{\perp}$
and $c^{\perp}$ are the following
\begin{align}
a^{\perp} & =\{(x:y:z:t) \mid x+(1-c)z-t=0 \},  \notag \\
b^{\perp} & =\{(x:y:z:t) \mid y+(1-a)x-t=0 \},  \notag \\
c^{\perp} & =\{(x:y:z:t) \mid z+(1-b)y-t=0 \}.  \notag
\end{align}
}

{The poles of these planes $P_{a},P_{b}$ and $P_{c}$ are given by
\begin{align}
P_{a} & =(-A^{2}:0:(c-1)C^{2}:1),  \notag \\ P_{b} &
=((a-1)A^{2}:-B^{2}:0:1),  \notag \\ P_{c} &
=(0:(b-1)B^{2}:-C^{2}:1)  \notag
\end{align}
respectively. Then the conditions of the orthogonality
\begin{equation*}
\cos(a^{\perp},b^{\perp})=\cos(b^{\perp},c^{\perp})=\cos(c^{\perp},a^{\perp
})=0,
\end{equation*}
are equivalent to the equalities
\begin{equation*}
(P_{a},P_{b})=(P_{b},P_{c})=(P_{c},P_{a})=0.
\end{equation*}
As a result we have the following relations
\begin{equation*}
(1-a)A^{2}+1=(1-b)B^{2}+1=(1-c)C^{2}+1=0
\end{equation*}
or
\begin{equation*}
a=1+\frac{1}{A^{2}},\ b=1+\frac{1}{B^{2}},\ c=1+\frac{1}{C^{2}}.
\end{equation*}
Hence
\begin{align}
P_{a}=(-A^{2}:0:1:1),  \notag \\
P_{b}=(1:-B^{2}:0:1),  \notag \\
P_{c}=(0:1:-C^{2}:1).  \notag
\end{align}
}

{The poles $a^{\prime},b^{\prime}$ and $c^{\prime}$ of the planes
$a^{\prime \perp}$, $b^{\prime \perp}$ and $c^{\prime \perp}$
symmetric to the planes $a^{\perp}$, $b^{\perp}$ and $c^{\perp}$
with respect to $z=0$, $x=0$ and $y=0$ respectively are given by
\begin{align}
P_{a^{\prime}}=(-A^{2}:0:-1:1),  \notag \\
P_{b^{\prime}}=(-1:-B^{2}:0:1),  \notag \\
P_{c^{\prime}}=(0:-1:-C^{2}:1).  \notag
\end{align}
}

{The angle between $a^{\perp}$ and ${a^{\prime}}^{\perp}$ is equal ${2
\alpha }$ and we have
\begin{equation*}
\cos{2\alpha}=-\frac{(P_{a},P_{a^{\prime}})} {\sqrt{{(P_{a},P_{a})(P_{a^{
\prime}},P_{a^{\prime}})}}} =-\frac{A^{2}-\frac{1}{C^{2}}+1}{A^{2}+\frac
{1}{C^{2}}+1}.
\end{equation*}
Hence
\begin{equation*}
\tan^{2}{\alpha}=C^{2}(A^{2}+1)
\end{equation*}
and by analogy we obtain
\begin{align}
\tan^{2}{\beta} & =A^{2}(B^{2}+1),  \notag \\
\tan^{2}{\gamma} & =B^{2}(C^{2}+1).  \notag
\end{align}
}

{As a result we have the following }

\begin{lemma}
\label{l1}{{{The polyhedron }$P(a,b,c)${ in the Klein model }}}$%
\mathbb{K}${{\ {of the spherical space }}}$\mathbb{S}${{$^{3}=(S^{3},{%
d\sigma }^{2})$ {has the right dihedral angles between the faces
intersecting at the point }}}$O^{\ast }${{{\ if and only if }
\begin{equation*}
a=1+\frac{1}{A^{2}},\ b=1+\frac{1}{B^{2}},\ c=1+\frac{1}{C^{2}}.
\end{equation*}
{Moreover, the dihedral angles }$\alpha$, $\beta${ and }$\gamma$ of $%
P(a,b,c)$ {satisfy the relations }
\begin{align}
\tan ^{2}{\alpha }=C^{2}(A^{2}+1),  \notag \\
\tan ^{2}{\beta }=A^{2}(B^{2}+1),  \notag \\
\tan ^{2}{\gamma }=B^{2}(C^{2}+1).  \notag
\end{align}
\vspace{3pt} }}
\end{lemma}

{{Set $L=\tan {\alpha }$, $M=\tan {\beta }$, $N=\tan {\gamma }$
and find an algebraic equation for variable
\begin{equation*}
T=-ABC
\end{equation*}
in terms of $L$, $M$, $N$. By Lemma \ref{l1}, we have
\begin{align}
C^{2}A^{2}+C^{2}& =L^{2},  \notag \\
A^{2}B^{2}+A^{2}& =M^{2},  \label{ff0} \\
B^{2}C^{2}+B^{2}& =N^{2}.  \notag
\end{align}
Multiplying the first of equations by $B^{2}$ and subtracting the
third, we express
$B^{2}$ by means of $A^{2}B^{2}C^{2}$. Similarly, expressing $A^{2}$ and $%
C^{2}$, we have
\begin{align}
A^{2}B^{2}C^{2}+N^{2}& =(1+L^{2})B^{2},  \notag \\
A^{2}B^{2}C^{2}+L^{2}& =(1+M^{2})C^{2},  \label{f0} \\
A^{2}B^{2}C^{2}+M^{2}& =(1+N^{2})A^{2}.  \notag
\end{align}
By multiplying the equations, we obtain
\begin{equation*}
(T^{2}+L^{2})(T^{2}+M^{2})(T^{2}+N^{2})=(1+L^{2})(1+M^{2})(1+N^{2})T^{2},
\end{equation*}
where $T=-ABC$. }}

{{Note that $T=\pm 1$ is a root of the above equation of order six. Rewrite
the equation in the following equivalent form
\begin{equation*}
(T^{2}-1)(T^{4}+(L^{2}+M^{2}+N^{2}+1)T^{2}-L^{2}M^{2}N^{2})=0.
\end{equation*}
In the case $T^{2}=1$ we have }}$A^{2}B^{2}C^{2}=1$. Then{{\ from{\ ({\ref
{f0}}) and (\ref{ff0}) we deduce }the identity
\begin{equation*}
1^{4}+(L^{2}+M^{2}+N^{2}+1)1^{2}-L^{2}M^{2}N^{2}=0.
\end{equation*}
Hence, the equation under consideration is equivalent to
\begin{equation*}
T^{4}+(L^{2}+M^{2}+N^{2}+1)T^{2}-L^{2}M^{2}N^{2}=0.
\end{equation*}
By definition, }}$T=-ABC$, where $A,B,C>0$. {{\ Hence, $T$ as a
negative root of the equation
\begin{equation*}
T^{2}=-\left( \frac{L^{2}+M^{2}+N^{2}+1}{2}\right) +\sqrt{{\left( \frac{%
L^{2}+M^{2}+N^{2}+1}{2}\right) }^{2}+L^{2}M^{2}N^{2}}.
\end{equation*}
}}

{{\ Notice that (\ref{f0}) can be rewritten in the form
\begin{equation*}
A^{2}=\frac{T^{2}+M^{2}}{1+N^{2}},\ B^{2}=\frac{T^{2}+N^{2}}{1+L^{2}},\
C^{2}=\frac{T^{2}+L^{2}}{1+M^{2}}.
\end{equation*}
As a result we obtain the following }}

\begin{lemma}
\label{l2}{{{The values }$A,B,C${\ and }$T=-ABC${\ satisfy the following
relations} {\ }
\begin{equation}
T^{4}+(L^{2}+M^{2}+N^{2}+1)T^{2}-L^{2}M^{2}N^{2}=0,  \notag
\end{equation}
\begin{equation*}
T^{2}=-p+\sqrt{p^{2}+L^{2}M^{2}N^{2}},
\end{equation*}
{where}
\begin{align}
p &=\frac{L^{2}+M^{2}+N^{2}+1}{2},  \notag \\
A^{2} &=
\frac{T^{2}+M^{2}}{1+N^{2}},\ B^{2}=\frac{T^{2}+N^{2}}{1+L^{2}},\
C^{2}=\frac{T^{2}+L^{2}}{1+M^{2}}.  \notag
\end{align}
{Here }$A,B,C>0${\ and }$T=-ABC<0${\ are uniquely determined by the values }$%
L=\tan \alpha ,M=\tan \beta ${\ and }$N=\tan \gamma ${. } }}
\end{lemma}

{{Now we shall find the lengths of edges for spherical polyhedron
$Q({\alpha },\beta ,\gamma )$. First we find the length
${L_{\alpha }}$ of the edge with dihedral angle ${\ \alpha }$. The
{projective coordinates {of the } }terminal points of the edge are
given by
\begin{equation*}
u_{100}=(1:0:0:1),\hspace{7mm}u_{1a0}=(1:a:0:1).
\end{equation*}
Let $s_{100}$ and $s_{1a0}$ be their projections on the upper hemisphere
\begin{equation*}
{S^{3}}_{+}=\{v=(x,y,z,t):(v,v)=1,t>0\}.
\end{equation*}
We have
\begin{equation*}
s_{100}=(\lambda :0:0:\lambda ),
\end{equation*}
where $\frac{1}{A^{2}}\lambda ^{2}+\lambda ^{2}=1$. Hence
\begin{equation*}
\lambda =\frac{A}{\sqrt{1+A^{2}}}.
\end{equation*}
By analogy, for
\begin{equation*}
s_{1a0}=(\tilde{\lambda}:\tilde{\lambda}a:0:\tilde{\lambda})\in {S^{3}}_{+}
\end{equation*}
we have
\begin{equation*}
\frac{\tilde{\lambda}^{2}}{A^{2}}+\frac{\tilde{\lambda}^{2}a^{2}}{B^{2}}+%
\tilde{\lambda}^{2}=1,
\end{equation*}
where $a=1+\frac{1}{A^{2}}$ is defined by Lemma \ref{l1}. It gives
\begin{equation*}
\tilde{\lambda}=\frac{AB}{\sqrt{A^{2}(B^{2}+1)+1}}\frac{A}{\sqrt{A^{2}+1}}.
\end{equation*}
}}

{{We obtain
\begin{equation*}
\cos{L_{\alpha}}=(s_{100},s_{1a0})=\frac{\tilde{\lambda}{\lambda}}{A^{2}}+%
\tilde{\lambda}{\lambda}=\frac{AB}{\sqrt{A^{2}(B^{2}+1)+1}}.
\end{equation*}
We notice that the polyhedron $P(a,b,c)$ is contained in the first
octant of the Euclidean space }}$\mathbb{R}${{$^{3}$(see Fig.
\ref{fig1}).
 Hence, the spherical length $L_{\alpha}$ of the edge whose dihedral angle $%
{\alpha}$ satisfies the inequality $0<L_{\alpha}<\frac{\pi}{2}$.
Thus
\begin{equation*}
\sin{L_{\alpha}}=\frac{\sqrt{A^{2}+1}}{\sqrt{A^{2}(B^{2}+1)+1}}
\end{equation*}
and
\begin{equation*}
\tan{L_{\alpha}}=\frac{\sqrt{A^{2}+1}}{AB}.
\end{equation*}
By analogy we obtain
\begin{align}
\tan{L_{\beta}} & =\frac{\sqrt{B^{2}+1}}{BC},  \notag \\
\tan{L_{\gamma}} & =\frac{\sqrt{C^{2}+1}}{AC},  \notag
\end{align}
So we prove the following }}

\begin{lemma}
\label{l3}{{\ {The spherical lengths }$L_{\alpha },L_{\beta }$ {\ and }$%
L_{\gamma }${\ are given by formulas }
\begin{equation*}
\tan {L_{\alpha }}=\frac{\sqrt{A^{2}+1}}{AB},\ \ \tan {L_{\beta }}=\frac{%
\sqrt{B^{2}+1}}{BC},\ \ \tan {L_{\gamma }}=\frac{\sqrt{C^{2}+1}}{AC}.\ \
\end{equation*}
}}
\end{lemma}

{As an immediate corollary of Lemma \ref{l1} and Lemma \ref{l2} we obtain
the following theorem (see also \cite{D}). }

\begin{theorem}
{{\ {Let }$\alpha ,\beta ${\ and} $\gamma ${\ are such that }$\frac{\pi }{2%
}<\alpha ,\beta ,\gamma <\pi ${. Then there exists a spherical Lambert cube
with the essential angles} $\alpha ,\beta $ {and }$\gamma $ {. } }}
\end{theorem}

 \begin{proof}
  { By Lemma \ref{l2}, the values
$A^{2},B^{2}$ and $C^{2}$ uniquely determined by the equalities
\begin{equation*}
A^{2}=\frac{T^{2}+M^{2}}{1+N^{2}},\ \ B^{2}=\frac{T^{2}+N^{2}}{1+L^{2}},\ \
C^{2}=\frac{T^{2}+L^{2}}{1+M^{2}}.
\end{equation*}
Then by Lemma \ref{l1}, we conclude the existence of the Euclidean
polyhedron $P(a,b,c)$. }

{{The spherical cube $Q(\alpha,\beta,\gamma)$ with the essential angles $%
\alpha,\beta$ and $\gamma$ is a result of projection of $P(a,b,c)$ into $%
S^{3}$. }}
 \end{proof}

{For our purposes we need metrical relations between essential angles and
lengths of Lambert cube $Q(\alpha,\beta,\gamma)$. They are given by the
following two theorems. }

\begin{theorem}[The Tangent Rule]
\label{t4}{{\ {Let {$Q(\alpha ,\beta ,\gamma )$,{\ {$%
\frac{\pi }{2}<\alpha ,\beta ,\gamma <\pi $} be {a spherical Lambert cube}}}}%
.${\ }$Denote by $L_{\alpha },L_{\beta }${\ and }$L_{\gamma }${\ lengths of
edges whose dihedral angles are }$\alpha ,\beta $ {and }$\gamma ${\
respectively.} {Then}
\begin{equation*}
\frac{\tan {\alpha }}{\tan {L_{\alpha }}}=\frac{\tan {\beta }}{\tan {%
L_{\beta }}}=\frac{\tan {\gamma }}{\tan {L_{\gamma }}}=T.
\end{equation*}
{where }$T$ {is a negative root of the equation}
\begin{equation*}
T^{4}+(L^{2}+M^{2}+N^{2}+1)T^{2}-L^{2}M^{2}N^{2}=0
\end{equation*}
{and }$L=\tan \alpha ,\ M=\tan \beta ,\ N=\tan \gamma $ . }}
\end{theorem}

  \begin{proof}
{{ By Lemma \ref{l1}, Lemma \ref{l3} and the condition $%
\frac{\pi}{2}<\alpha,\beta,\gamma<\pi$, we have
\begin{equation*}
\tan{\alpha}=-C\sqrt{A^{2}+1}
\end{equation*}
and
\begin{equation*}
\tan{L_{\alpha}}=\frac{\sqrt{A^{2}+1}}{AB}.
\end{equation*}
Then
\begin{equation*}
\frac{\tan{\alpha}}{\tan{L_{\alpha}}}=-ABC=T,
\end{equation*}
where $T$ is the same as in Lemma \ref{l2}. By the same way, we
obtain the other equalities of the theorem. }}
 \end{proof}

{{The parameter $T$ from Theorem \ref{t4} can be presented in the form $%
T=\tan{\theta}$, where $\frac{\pi}{2}<\theta<\pi$. The value $\theta$ plays
the significant role in the studying of metric structure of Lambert cube and
is called a \textit{principal parameter}. }}

{The following theorem is a direct corollary of Theorem \ref{t4}. }

\begin{theorem}[The Sine-Cosine Rule]
\label{t5}{\ { {{Let {$Q(\alpha ,\beta ,\gamma )$,{\ {$\frac{\pi
}{2}<\alpha ,\beta ,\gamma <\pi $} be {a
spherical Lambert cube}}}}.${\ }$Denote by $L_{\alpha },L_{\beta }${\ and }$%
L_{\gamma }${\ lengths of edges whose dihedral angles }$\alpha
,\beta $ {and }$\gamma ${,  respectively.} Then }
\begin{equation*}
\frac{\sin {\alpha }}{\sin {L_{\alpha }}}{\cdot }\frac{\sin {\beta }}{\sin {%
L_{\beta }}}{\cdot }\frac{\cos {\gamma }}{\cos {L_{\gamma }}}=-1.
\end{equation*}
}}
\end{theorem}

 \begin{proof}
{ From Theorem \ref{t4}, we have
\begin{equation*}
\tan {L_{\alpha }}=\frac{L}{T},\ \ \tan {L_{\beta }}=\frac{M}{T},\ \ \tan {\
L_{\gamma }}=\frac{N}{T}.
\end{equation*}
Hence
\begin{equation*}
\sin ^{2}{L_{\alpha }}=\frac{L^{2}}{T^{2}+L^{2}},\ \ \sin ^{2}{L_{\beta }}=%
\frac{M^{2}}{T^{2}+M^{2}},\ \ \cos ^{2}{L_{\gamma }}=\frac{T^{2}}{T^{2}+N^{2}%
}
\end{equation*}
and by Lemma \ref{l2}, we have
\begin{equation}
\frac{\sin ^{2}{\alpha }}{\sin ^{2}{L_{\alpha }}}{\cdot }\frac{\sin ^{2}{%
\beta }}{\sin ^{2}{L_{\beta }}}{\cdot }\frac{\cos ^{2}{\gamma }}{\cos ^{2}{%
L_{\gamma }}}=\frac{(T^{2}+L^{2})(T^{2}+M^{2})(T^{2}+N^{2})}{%
(1+L^{2})(1+M^{2})(1+N^{2})T^{2}}=1.  \label{scr}
\end{equation}
Since $\frac{\pi }{2}<\alpha ,\beta ,\gamma <\pi $ and the values $%
0<L_{\alpha },L_{\beta },L_{\gamma }<\frac{\pi }{2}$, the theorem is proved. }%

 \end{proof}

\begin{remark}
{{\ U{p to cyclic permutation of angles }}}$\alpha ,\beta $ and $\gamma $,{{{%
\ }Theorem \ref{t5} contains three independent equations which are
sufficient to determine $L_{\alpha }$, $L_{\beta }$ and $L_{\gamma }$ in
terms of ${\alpha }$, ${\beta }$ and ${\gamma }$. }}
\end{remark}

\section{{{{The volume of spherical Lambert cube} }}}

{{From now on, for an arbitrary function $\phi:(\frac{\pi}{2},\pi ){%
\rightarrow}$}}$\mathcal{R}${{\ we shall write $\phi(\frac{\pi}{2})$ and $%
\phi({\pi})$ instead of $\phi(\frac{\pi}{2}+0)$ and $\phi({\pi}-0)$,
respectively, if the corresponding values have a sense. }}

{Let $V=V(\alpha,\beta,\gamma)$ be the volume of spherical Lambert cube $%
Q(\alpha,\beta,\gamma)$ with essential angles $\alpha,\beta$ and
$\gamma$. }

{{Notice first that $V(\alpha, \beta, \gamma) \rightarrow0$, as $\alpha,
\beta, \gamma\rightarrow\frac{\pi}{2}+0$. }}

{{It can be justified by the following arguments. Each face of $Q(\alpha
,\beta ,\gamma )$ is a Lambert quadrilateral with essential angles ${\alpha }
$, ${\beta }$ and ${\gamma }$ tending to $\frac{\pi }{2}$. Then there is {a
pair of opposite sides{\ in the quadrilateral whose}} lengths tend to zero.
We have decrease in the dimension in each face of $Q(\alpha ,\beta ,\gamma )$
and hence of $Q(\alpha ,\beta ,\gamma )$ itself. It means the volume tends
to zero. }}

{{Thus, we shall consider the equality
\begin{equation}
V(\frac{\pi}{2},\frac{\pi}{2},\frac{\pi}{2})=0  \label{v0}
\end{equation}
as a true one. }}

{{Let ${L_{\alpha }},{L_{\beta }}$ and ${L_{\gamma }}$ be the
lengths of edges whose dihedral angles ${\alpha },\beta $ and
$\gamma $, respectively. By the Schl\"{a}fli formula (see
\cite{Sh}, \cite{M}), we have
\begin{equation}
\frac{\partial V}{\partial \alpha }=\frac{1}{2}{\ }L_{\alpha }{\,},\qquad \frac{%
\partial V}{\partial \beta }=\frac{1}{2}{\ }L_{\beta }{\,},\qquad \frac{\partial
V}{\partial \gamma }=\frac{1}{2}{\ }L_{\gamma }{\,}.  \label{f2}
\end{equation}
}}

\begin{proposition}
\label{main} \ The volume of spherical Lambert cube $%
Q(\alpha ,\beta ,\gamma ),${\ }$\frac{\pi }{2}<\alpha ,\beta ,\gamma <\pi $
is given by the formula
\begin{equation}
V(\alpha ,\beta ,\gamma )={\frac{1}{4}\int\limits_{-\infty }^{T}{\log {\frac{%
(t^{2}+L^{2})(t^{2}+M^{2})(t^{2}+N^{2})}{(1+L^{2})(1+M^{2})(1+N^{2})t^{2}}}%
\frac{dt}{t^{2}-1},}}  \label{v}
\end{equation}
where $T$ is a negative root of the equation
\begin{equation*}
T^{4}+(L^{2}+M^{2}+N^{2}+1)T^{2}-L^{2}M^{2}N^{2}=0,
\end{equation*}
$L=\tan \alpha ,\ M=\tan \beta \ $and $N=\tan \gamma .$
\end{proposition}

 \begin{proof}
 By the above arguments the volume function $V=V(\alpha ,\beta ,\gamma
)$ satisfies $(\ref{f2})$ with initial data $(\ref{v0})$. We set
\begin{equation*}
\widetilde{V}={\int\limits_{-\infty }^{T}F(t,L,M,N) \ dt{,}}
\end{equation*}
where
\begin{equation*}
F(t,L,M,N)={\frac{1}{4(t^{2}-1)}{\log {\frac{%
(t^{2}+L^{2})(t^{2}+M^{2})(t^{2}+N^{2})}{(1+L^{2})(1+M^{2})(1+N^{2})t^{2}}}}}
\end{equation*}
and show that \ $\widetilde{V}$ satisfies conditions $(\ref{v0})$ and $(\ref
{f2})$ . Then $V=\widetilde{V}$ and the proposition is proven.

By the Leibnitz formula we obtain
\begin{equation}
\frac{\partial \widetilde{{V}}}{\partial {\alpha }}=F(T,L,M,N)\frac{\partial
T}{\partial {\alpha }}+{\int\limits_{-\infty }^{T}\frac{\partial F(t,L,M,N)}{%
\partial {L}}\frac{\partial L}{\partial {\alpha }}dt}  \label{f4}
\end{equation}

By $\left( \ref{scr}\right) $ we have $F(T,L,M,N)=0$. Moreover, since ${%
\alpha =\arctan }${$L$}, we have
\begin{equation*}
{\frac{\partial L}{\partial {\alpha }}=1+L^{2}}\quad \text{and\quad }{\frac{%
\partial F(t,L,M,N)}{\partial {L}}\frac{\partial L}{\partial {\alpha }}=-}%
\frac{L}{2\left( t^{2}+L^{2}\right) }\text{.}
\end{equation*}
Hence, by Tangent Rule, we obtain from $(\ref{f4})$
\begin{equation*}
\frac{\partial \widetilde{{V}}}{\partial {\alpha }}={\int\limits_{-\infty
}^{T}{-}\frac{L{dt}}{2\left( t^{2}+L^{2}\right) }=}\frac{1}{2} \ {\arctan }%
\frac{L}{T}=\frac{1}{2} \ L_{\alpha }.
\end{equation*}
The equalities
\begin{equation*}
\frac{\partial \widetilde{{V}}}{\partial {\beta }}=\frac{1}{2} \
L_{\beta
}\quad \text{and\quad }\frac{\partial \widetilde{{V}}}{\partial {\gamma }}=%
\frac{1}{2} \ L_{\gamma }
\end{equation*}
can be obtained by the similar way.

To verify the initial condition $(\ref{v0})$ for function
$\widetilde{{V}}$ we remark that $L,M,N\rightarrow +\infty $ as
{{$\alpha ,\beta ,\gamma \rightarrow \frac{\pi }{2}+0$, hence
}}$T\rightarrow -\infty $. From the convergence of integral
$(\ref{v})$, we have $\widetilde{{V}}\rightarrow 0$ as
$T\rightarrow -\infty $.
\end{proof}

\section{The proof of the main theorem}

\begin{proof}[Proof of Theorem \ref{mm}]
By Proposition \ref{main}, we have
\begin{align}
V\left( \alpha ,\beta ,\gamma \right) & ={}\frac{1}{4}\int\limits_{-\infty
}^{T}{\log }\left( \frac{t^{2}+L^{2}}{1+L^{2}}\frac{t^{2}+M^{2}}{1+M^{2}}{{%
\frac{t^{2}+N^{2}}{1+N^{2}}:}}\frac{t^{2}+0^{2}}{1+0^{2}}\right) {\frac{dt}{%
t^{2}-1}}  \label{v1} \\
& =\frac{1}{4}\left( I\left( L,T\right) +I\left( M,T\right) +I\left(
N,T\right) -I\left( 0,T\right) \right)  \notag
\end{align}
where
\begin{equation*}
I\left( L,T\right) =\int\limits_{-\infty }^{T}{\log }\left( {{\frac{%
t^{2}+L^{2}}{1+L^{2}}}}\right) {\frac{dt}{t^{2}-1}.}
\end{equation*}
Let $T=\tan \theta $, $L=\tan \alpha $, $M=\tan \beta $ and $N=\tan \gamma $%
. Then under substitution $t=\tan \tau $ we obtain
\begin{align*}
I\left( L,T\right) & =\int\limits_{\frac{\pi }{2}}^{\theta }{\log }\left( {{%
\frac{\tan ^{2}\tau +\tan ^{2}\alpha }{1+\tan ^{2}\alpha }}}\right) {\frac{%
d\tau }{\cos ^{2}\tau \left( \tan ^{2}\tau -1\right) }} \\
{}& =\int\limits_{\theta }^{\frac{\pi }{2}}{\frac{{\log }\left( 1-\cos 2\tau
\cos 2\alpha \right) d\tau }{\cos 2\tau }}-\int\limits_{\theta }^{\frac{\pi
}{2}}{\frac{{\log }\left( 1+\cos 2\tau \right) d\tau }{\cos 2\tau }} \\
& =\delta \left( \alpha ,\theta \right) -\delta \left( \frac{\pi }{2},\theta
\right) .
\end{align*}
Hence, (\ref{v1}) yields
\begin{align*}
V\left( \alpha ,\beta ,\gamma \right) & =\frac{1}{4}\left( I\left(
L,T\right) +I\left( M,T\right) +I\left( N,T\right) -I\left( 0,T\right)
\right) \\
& =\frac{1}{4}\left( \delta \left( \alpha ,\theta \right) +\delta \left(
\beta ,\theta \right) +\delta \left( \gamma ,\theta \right) -2\delta \left(
\frac{\pi }{2},\theta \right) -\delta \left( {0},\theta \right) \right) .
\end{align*}

\end{proof}

\section{Explicit volume calculations}

{This section is devoted to the explicit volume calculations in some
particular cases. }

\begin{proposition}
\label{p3}{{{Let } $\alpha ,\beta ,\gamma ${,} $\frac{\pi
}{2}<\alpha,\beta,\gamma<\pi $ {are related by the following
equation}
\begin{equation*}
\cos ^{2}\alpha +\cos \beta ^{2}+\cos \gamma ^{2}={1}\text{.}
\end{equation*}
{Then the volume of a spherical Lambert cube} $Q(\alpha ,\beta ,\gamma )$ {\
is given by the formula }
\begin{equation*}
V(\alpha ,\beta ,\gamma )=\frac{1}{4}\left( \frac{{\pi }^{2}}{2}-{(\pi
-\alpha )}^{2}-{(\pi -\beta )}^{2}-{(\pi -\gamma )}^{2}\right) .
\end{equation*}
}}
\end{proposition}

\begin{proof}
Remind that
\begin{equation*}
\tan \theta =-\sqrt{-p+\sqrt{p^{2}+L^{2}M^{2}N^{2}}},
\end{equation*}
{where}
\begin{equation*}
{\ }p=\frac{L^{2}+M^{2}+N^{2}+1}{2}
\end{equation*}
{{{and }$L=\tan \alpha ,\ M=\tan \beta ,\ N=\tan \gamma $ . }Let }
\begin{equation*}
\cos ^{2}\alpha +\cos \beta ^{2}+\cos \gamma ^{2}={1}\text{.}
\end{equation*}
By applying the elementary trigonometry we have
\begin{equation*}
2p+1=\frac{1}{\cos ^{2}\alpha }+\frac{1}{\cos ^{2}\beta }+\frac{1}{\cos
^{2}\gamma }-1
\end{equation*}
and
\begin{eqnarray*}
L^{2}M^{2}N^{2} &=&\frac{1}{\cos ^{2}\alpha }+\frac{1}{\cos ^{2}\beta }+%
\frac{1}{\cos ^{2}\gamma }-1 \\
&&-\frac{\cos ^{2}\alpha +\cos ^{2}\beta +\cos ^{2}\gamma -1}{\cos
^{2}\alpha \cos ^{2}\beta \cos ^{2}\gamma } \\
&=&2p+1\text{.}
\end{eqnarray*}
Hence{{\ $\tan \theta =-1$ and consequently, $\theta =\frac{3\pi
}{4}$. By Theorem \ref{mm} and Corollary \ref{td1}}}${(iii)}$\
(see Appendix), we obtain
\begin{align*}
V\left( \alpha ,\beta ,\gamma \right) & =\frac{1}{4}\left( \delta \left(
\alpha ,\theta \right) +\delta \left( \beta ,\theta \right) +\delta \left(
\gamma ,\theta \right) -2\delta \left( \frac{\pi }{2},\theta \right) -\delta
\left( {0},\theta \right) \right) \\
{}& =\frac{1}{4}\left( \frac{{\pi }^{2}}{2}-{(\pi -\alpha )}^{2}-{(\pi
-\beta )}^{2}-{(\pi -\gamma )}^{2}\right) .
\end{align*}
\end{proof}

In particular, we have the following.

\begin{corollary}
\begin{equation*}
V\left( \frac{2\pi }{3},\frac{2\pi }{3},\frac{3\pi }{4}\right) =\frac{31}{576%
}\pi ^{2}.
\end{equation*}
\end{corollary}

{{Consider the process of degenerating of Lambert cube $Q(\alpha
,\beta ,\gamma )$ as $\gamma \rightarrow \pi $. Let $L_{\gamma
}^{\ast }$ be the length of the edge opposite to the edge with
essential angle $\gamma $. The straightforward calculations based
on Theorem \ref{t5} show that $L_{\gamma }^{\ast }$ tends to zero
as $\gamma \rightarrow \pi $. As a result the {singular } Lambert
cube $Q(\alpha ,\beta ,\pi )$ can be defined as a cone under the
plane spherical quadrilateral with the angles $\frac{\pi }{2},\alpha ,\frac{%
\pi }{2},\beta $ (Fig \ref{fig2}) . }

 \begin{figure}
 \resizebox{1.05\hsize}{!}{\includegraphics*{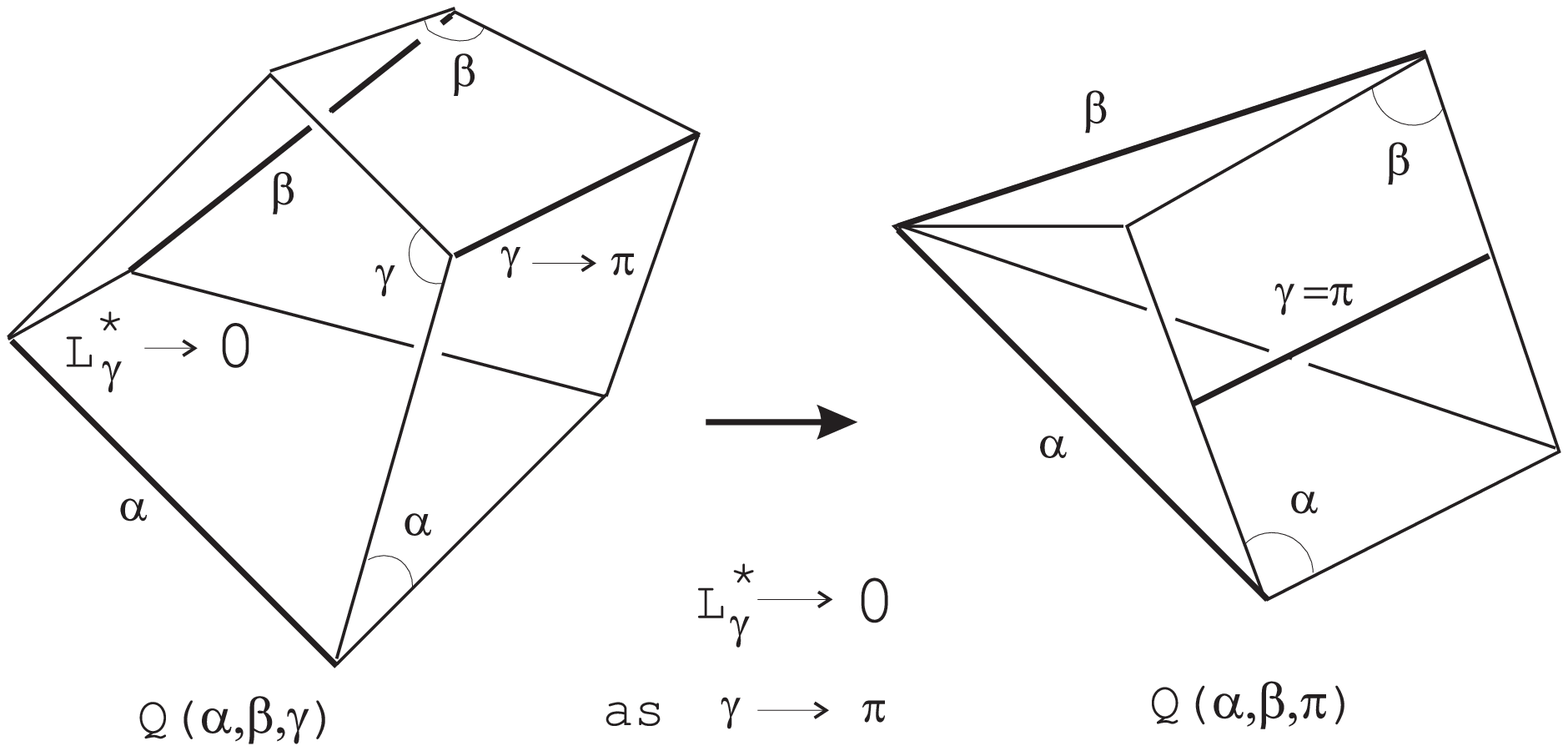}}
 \caption{} \label{fig2}
 \end{figure}

{Applying Theorem \ref{mm} and Corollary \ref{td1}}${(iv)}${\ in case $%
\gamma =\pi $, we obtain the following }

\begin{proposition}
{{{The volume of a singular Lambert cube} $Q(\alpha ,\beta ,\pi )$, $\frac{%
\pi }{2}<\alpha ,\beta <\pi $ {is given by the formula}
\begin{equation*}
V(\alpha ,\beta ,\pi )=\left( \alpha +\beta -\pi \right) {\pi }.
\end{equation*}
}}
\end{proposition}

\section{Appendix}

\bigskip This section is devoted to elementary properties of the function {$%
\delta \left( \alpha ,\theta \right) $ and its relations to the
Dilogarithm, Lobachevskij and Schl\"{a}fli functions.}

\subsection{\protect\bigskip {Elementary properties of $\protect\delta %
\left( \protect\alpha ,\protect\theta \right) $}}

{\ }Now we list the {following} elementary properties of the
function ${\delta \left( \alpha ,\theta \right) }$.

\begin{proposition}
The function
\begin{equation*}
{\delta \left( \alpha ,\theta \right) }=\int\limits_{\theta }^{\frac{\pi }{2}%
}{\log }\left( {1-\cos 2\alpha \cos 2\tau }\right) \frac{d\tau }{\cos 2\tau }
\end{equation*}
satisfies the following properties

\begin{tabular}{ll}
$(i)$ & {$\delta \left( \alpha ,\theta \right) $ is continuous for
all }$(\alpha ,\theta )\in R^{2}$ and \\ & differentiable with
respect to{\ }$(\alpha ,\theta )$ for $\alpha \neq \pi /2+k\pi
,k\in Z$; \\
$(ii)$ & {$\delta \left( \alpha ,\theta \right) $ is even} with respect to $%
\alpha $ and satisfies the relation \\ & $\delta \left( \alpha
,\theta \right) +\delta \left( \alpha ,-\theta \right) =2\delta
\left( \alpha ,0\right) $ for all{\ }$(\alpha ,\theta )\in R^{2}$;
\\ $(iii)$ & {$\delta \left( \alpha ,\theta \right) =$ $\delta
\left( \pi
-\alpha ,\theta \right) $ }and {$\delta \left( \alpha ,\theta \right) =-$ $%
\delta \left( \alpha ,\pi -\theta \right) $}; \\ $(iv)$ & {$\delta
\left( \alpha ,\theta \right) $ is }$\pi $-periodic with respect
to $\alpha $ ; \\ $(v)$ & {$\delta \left( \alpha ,\theta \right) $
is }linear periodic with respect to $\theta $ in the following
\\ & sense \\ & $\delta \left( \alpha ,\theta +k\pi \right)
=\delta \left( \alpha ,\theta \right) -2k\delta \left( \alpha
,0\right) ,k\in Z$ \\ & and \\ & {{$\delta \left( \alpha ,0\right)
=\pi ^{2}/4-\left| \pi ^{2}/2-\alpha \pi \right| ,0\leq $}}$\alpha
\leq \pi $; \\ $(vi)$ & Let $\widetilde{\delta }\left( \alpha
,\theta \right) =\delta \left( \alpha ,\theta \right) +\left(
2\theta /\pi -1\right) \delta \left( \alpha ,0\right) $. Then \\
& $a)\quad \widetilde{\delta }\left( \alpha ,\theta \right) $ {is even} and{%
\ }$\pi $-periodic with respect to $\alpha $; \\ & $b)\quad
\widetilde{\delta }\left( \alpha ,\theta \right) $ {is odd} and{\
}$\pi $-periodic with respect to $\theta $; \\ & $c)\quad \left|
\widetilde{\delta }\left( \alpha ,\theta \right) \right| \leq \pi
^{2}/4$ and $\widetilde{\delta }\left( \pi /2,3\pi /4\right) =\pi
^{2}/4$.
\end{tabular}
\end{proposition}

\bigskip The {properties }$(i)-(iv)${\ can be deduced directly from the
definition of } {$\delta \left( \alpha ,\theta \right) $. To prove
the second statement of the property }${(ii)}$ we observe that the
equality is evidently true for $\theta =0$.{\ By definition of
$\delta \left( \alpha ,\theta \right) $, the derivatives of the
both sides of the equality with respect to }$\theta ${\ are equal
to zero.\ To see the equality}
\begin{equation*}
\delta \left( \alpha ,\theta +k\pi \right) =\delta \left( \alpha ,\theta
\right) -2k\delta \left( \alpha ,0\right) ,k\in Z
\end{equation*}
we notice that the derivatives of the both sides with respect to
$\theta $ are equal. To establish the property we have to check
the equality for one fixed value of $\theta $. Let $\theta =0$.
The equality holds for $k=0$. For any $k\in Z$ from properties
$(iii)$ and $(ii)$ we have
\begin{equation*}
\delta \left( \alpha ,\left( k+1\right) \pi \right) =-\delta \left( \alpha
,-k\pi \right)
\end{equation*}
and
\begin{equation*}
-\delta \left( \alpha ,-k\pi \right) =\delta \left( \alpha ,k\pi \right)
-2\delta \left( \alpha ,0\right) \text{.}
\end{equation*}
Hence
\begin{equation*}
\delta \left( \alpha ,\left( k+1\right) \pi \right) =\delta \left( \alpha
,k\pi \right) -2\delta \left( \alpha ,0\right)
\end{equation*}
and the first statement of $(v)$ follows by induction. The value of ${\delta
\left( \alpha ,0\right) }$ can be taken from the corollary \ref{td1} bellow.
The property $(vi)$ follows from properties $(i)-(iv)$.

For explicit calculations it is more convenient to use the following
alternative form of {$\delta \left( \alpha ,\theta \right) $. }

\begin{proposition}
\label{p2}{{{Let } {${\pi /2}<\alpha ,\theta <\pi $, then the function } $%
\delta \left( \alpha ,\theta \right) $ can be represented in the form }}
\begin{equation*}
{{\delta \left( \alpha ,\theta \right) =2\int\limits_{\frac{3\pi }{4}%
}^{\alpha }arccot\hspace{1mm}\left( \frac{\cot {\nu }}{\cot \theta }\right)
d\nu ,}}
\end{equation*}
{{where $0<arccot{\ }x<{\pi }$ for all $x$. }}
\end{proposition}

\begin{proof}
{ By definition, we have
\begin{equation*}
\delta \left( \alpha ,\theta \right) =\int\limits_{\theta }^{\frac{\pi }{2}}{%
\log }\left( {1-\cos 2\alpha \cos 2\tau }\right) \frac{d\tau }{\cos 2\tau }.
\end{equation*}
Hence
\begin{align*}
\frac{\partial \delta \left( \alpha ,\theta \right) }{\partial \alpha }&
=\int\limits_{\theta }^{\frac{\pi }{2}}\frac{2{\sin 2\alpha \cos 2\tau }}{{%
1-\cos 2\alpha \cos 2\tau }}\frac{d\tau }{\cos 2\tau } \\
& =\pi -2\arctan \left( \cot {\alpha }\tan \theta \right) =2\,arccot{\ }%
\left( \cot {\alpha }\tan \theta \right) .
\end{align*}
Notice that
\begin{equation*}
\delta \left( \frac{3\pi }{4},\theta \right) =0.
\end{equation*}
Hence
\begin{equation*}
\delta \left( \alpha ,\theta \right) =\int\limits_{\frac{3\pi }{4}}^{\alpha }%
\frac{\partial \delta \left( \nu ,\theta \right) }{\partial \nu }d\nu
=2\int\limits_{\frac{3\pi }{4}}^{\alpha }arccot\left( \frac{\cot {\nu }}{%
\cot \theta }\right) d\nu .
\end{equation*}
}
\end{proof}

We use the proposition to calculate the value of $\delta \left(
\alpha ,\theta \right) $ in some cases. The results of our calculations are
collected in the following corollary.

\begin{corollary}
\label{td1}{{{{For any $\alpha $, $0\leq \alpha \leq {\pi }$ we have}}}}

\begin{tabular}{ll}
{$(i)$} & {{$\delta \left( \alpha ,0\right) =\pi \left( \frac{\pi }{4}%
-\left| \frac{\pi }{2}-\alpha \right| \right) $}}$,$ \\
{$(ii)$} & $\delta \left( \alpha ,\frac{\pi }{4}\right) =\left( \frac{\pi }{2%
}-\left| \frac{\pi }{2}-\alpha \right| \right) ^{2}-\frac{\pi ^{2}}{16},$ \\
{$(iii)$} & $\delta \left( \alpha ,\frac{3\pi }{4}\right) =\frac{\pi ^{2}}{16%
}-\left( \frac{\pi }{2}-\left| \frac{\pi }{2}-\alpha \right| \right) ^{2}${$%
, $} \\
$(iv)$ & $\delta \left( \alpha ,\pi \right) =\pi \left( \left| \alpha -\frac{%
\pi }{2}\right| -\frac{\pi }{4}\right) $.
\end{tabular}
\end{corollary}

{\smallskip The result easy follows from Proposition \ref{p2} and
properties }$(i)-(iii)$ of function {$\delta $.}

\subsection{The relation $\protect\delta \left( \protect\alpha ,\protect%
\theta \right) $ to the {Schl\"{a}fli }function}

\bigskip \bigskip Let $T(\alpha ,\beta ,\gamma )$ be a double-rectangular
spherical, Euclidean or hyperbolic tetrahedron with dihedral angles $\pi
/2-\alpha $, $\beta $ and $\pi /2-\gamma $. Define the {Schl\"{a}fli
function by the following formula}
\begin{equation*}
S(\alpha ,\beta ,\gamma )=\sum_{1}^{\infty }\frac{\left( -X\right) ^{n}}{%
n^{2}}\left( \cos 2n\alpha -\cos 2n\beta +\cos 2n\gamma -1\right) -\alpha
^{2}+\beta ^{2}-\gamma ^{2}\text{,}
\end{equation*}
where
\begin{equation*}
X=\frac{\sin \alpha \sin \gamma -D}{\sin \alpha \sin \gamma +D}\text{,}
\end{equation*}
\begin{equation*}
D=\sqrt{\cos ^{2}\alpha \cos ^{2}\gamma -\cos ^{2}\beta }\text{,}
\end{equation*}
\begin{equation*}
0\leq \alpha \leq \frac{\pi }{2},\qquad 0\leq \beta \leq \pi ,\qquad 0\leq
\gamma \leq \frac{\pi }{2}
\end{equation*}
(see Coxeter $[C]$ for details). It was shown by {Schl\"{a}fli (}$1898${)
that }in the spherical case ($\cos ^{2}\alpha \cos ^{2}\gamma >\cos
^{2}\beta $) the volume of $T(\alpha ,\beta ,\gamma )$ and $S(\alpha ,\beta
,\gamma )$ are related by {\ }
\begin{equation*}
4Vol(T(\alpha ,\beta ,\gamma ))=S(\alpha ,\beta ,\gamma )\text{.}
\end{equation*}

In Euclidean case ($\cos ^{2}\alpha \cos ^{2}\gamma =\cos ^{2}\beta $, see $%
[C]$, $p.~16$)
\begin{equation*}
S(\alpha ,\beta ,\gamma )=0\text{.}
\end{equation*}

 In hyperbolic case ($\cos ^{2}\alpha \cos ^{2}\gamma <\cos
^{2}\beta ,$ $\alpha ,\gamma <\beta $), Coxeter notice ( $[C]$ $,$
$p.~27$), that
\begin{equation*}
iS(\alpha ,\beta ,\gamma )=4Vol(T(\alpha ,\beta ,\gamma )).
\end{equation*}
From the other side, we know (see, for instance $[C]$, $p.~23$ or $[V]$, $%
p.~125$) that
\begin{eqnarray*}
4Vol(T(\alpha ,\beta ,\gamma )) &=&\Lambda (\alpha +\theta )+\Lambda
(-\alpha +\theta )-\Lambda (\beta +\theta )-\Lambda (-\beta +\theta ) \\
&&+\Lambda (\gamma +\theta )+\Lambda (-\gamma +\theta )+2\Lambda (\theta ),
\end{eqnarray*}
where $\theta $ is defined by
\begin{equation*}
\tan \theta =\frac{i\sin \alpha \sin \gamma }{D}=\frac{\sin \alpha \sin
\gamma }{\sqrt{\cos ^{2}\beta -\cos ^{2}\alpha \cos ^{2}\gamma }}\text{.}
\end{equation*}
Hence, the function $S(\alpha ,\beta ,\gamma )$ is connected with
Lobachevskij function by the relation
\begin{equation*}
iS(\alpha ,\beta ,\gamma )=-\Delta (\alpha ,\theta )+\Delta (\beta ,\theta
)-\Delta (\gamma ,\theta )+\Delta (0,\theta )\text{,}
\end{equation*}
where $\Delta (\alpha ,\theta )=\Lambda (\alpha +\theta )-\Lambda (\alpha
-\theta )$. The following proposition gives relation $S(\alpha ,\beta
,\gamma )$ and $\delta \left( \alpha ,\theta \right) $ in spherical case.

\begin{proposition}
Let $0\leq \alpha ,\beta ,\gamma \leq \frac{\pi }{2}$ and
\begin{equation*}
\tan \theta =\frac{\sin \alpha \sin \gamma }{\sqrt{\cos ^{2}\alpha \cos
^{2}\gamma -\cos ^{2}\beta }}\text{,}
\end{equation*}
where $\cos ^{2}\alpha \cos ^{2}\gamma >\cos ^{2}\beta $. Then
\begin{equation*}
S(\alpha ,\beta ,\gamma )=-\delta (\alpha ,\theta )+\delta (\beta ,\theta
)-\delta (\gamma ,\theta )+\delta (0,\theta )\text{.}
\end{equation*}
\end{proposition}

\begin{proof}
 Let $T(\alpha ,\beta ,\gamma )$ be a double-rectangular spherical
tetrahedron with essential dihedral angles $\frac{\pi }{2}-\alpha $, $\beta $
and $\frac{\pi }{2}-\gamma $. By the {Schl\"{a}fli theorem }$S(\alpha ,\beta
,\gamma )=4Vol\left( T(\alpha ,\beta ,\gamma )\right) $. We find $Vol\left(
T(\alpha ,\beta ,\gamma )\right) $ in a few steps.

{\it 1 step}. The Tangent Rule (see, for instance $[V]$, $p.125$)
\begin{equation*}
\frac{\tan \alpha }{\tan a}=\frac{\tan \beta }{\tan b}=\frac{\tan \gamma }{%
\tan c}=T\text{,}
\end{equation*}
where $a,b$ and $c$ are the spherical lengths of edges correspondent to
dihedral angles $\frac{\pi }{2}-\alpha ,\beta $ and $\frac{\pi }{2}-\gamma $
respectively,
\begin{equation*}
T=\frac{\sin \alpha \sin \gamma }{D}
\end{equation*}
and $D=\sqrt{\cos ^{2}\alpha \cos ^{2}\gamma -\cos ^{2}\beta }$.

{\it 2 step}. We note that $T$ is a root of the following
biquadratic equation
\begin{equation}
\frac{1+A^{2}}{T^{2}+A^{2}}\frac{T^{2}+B^{2}}{1+B^{2}}\frac{1+C^{2}}{%
T^{2}+C^{2}}T^{2}=1\text{,}  \label{F}
\end{equation}
where $A=\tan \alpha $, $B=\tan \beta $ and $C=\tan \gamma $.
Indeed, the above equation has four roots $T_{1,2}=\pm 1$ and
$T_{3,4}=\pm \sin \alpha \sin \gamma /D$. The equation $\left(
\ref{F}\right) $ has the following geometrical sense. By the
Tangent Rule, we have $\cos ^{2}a=T^{2}/\left(
T^{2}+A^{2}\right) $, $\cos ^{2}b=T^{2}/\left( T^{2}+B^{2}\right) $ and $%
\cos ^{2}c=T^{2}/\left( T^{2}+C^{2}\right) $. Hence $\left( \ref{F}\right) $
gives
\begin{equation*}
 \frac{\cos ^{2}a}{\cos ^{2}\alpha }{\cdot}\frac{\cos ^{2}\beta }{\cos ^{2}b}%
\cdot \frac{\cos ^{2}c}{\cos ^{2}\gamma }=1
\end{equation*}
and after the suitable choosing the sign we obtain
\begin{equation}
\frac{\cos \beta }{\cos b}=\frac{\cos \alpha }{\cos a}\cdot \frac{\cos
\gamma }{\cos c}\text{.}  \label{sc}
\end{equation}
This is the Cosine Rule for double-rectangular spherical tetrahedron $%
T(\alpha ,\beta ,\gamma )$.

{\it 3 step}. Let $V=V(\alpha ,\beta ,\gamma )$ be the volume of
$T(\alpha ,\beta ,\gamma )$. By the {{Schl\"{a}fli formula we
have}}
\begin{equation}
\frac{\partial V}{\partial {\alpha }}=-\frac{a}{2},\quad \frac{\partial V}{%
\partial {\beta }}=\frac{b}{2},\quad \frac{\partial V}{\partial {\gamma }}=-%
\frac{c}{2}\text{.}  \label{vt}
\end{equation}
It follows from the Tangent Rule that
\begin{equation}
V\rightarrow 0\quad \text{as\quad }T\rightarrow +\infty \text{.}  \label{it}
\end{equation}
{\it 4 step}. It is easy to check (see the proof of Proposition
{{\ref{main}}}) that the function
\begin{equation*}
\widetilde{V}(\alpha ,\beta ,\gamma )={\frac{1}{4}\int\limits_{T}^{+\infty }{%
\log {\frac{(1+A^{2})(t^{2}+B^{2})(1+C^{2})t^{2}}{%
(t^{2}+A^{2})(1+B^{2})(t^{2}+C^{2})}}\frac{dt}{t^{2}-1}}}
\end{equation*}
is a solution of $(\ref{vt})$ with the initial condition $(\ref{it})$. Hence
$Vol\left( T(\alpha ,\beta ,\gamma )\right) =\widetilde{V}(\alpha ,\beta
,\gamma )$. Since $S(\alpha ,\beta ,\gamma )=4Vol\left( T(\alpha ,\beta
,\gamma )\right) $ we have
\begin{eqnarray*}
{\ }S(\alpha ,\beta ,\gamma ) &=&{\int\limits_{T}^{+\infty }{\log {\frac{%
(1+A^{2})(t^{2}+B^{2})(1+C^{2})t^{2}}{(t^{2}+A^{2})(1+B^{2})(t^{2}+C^{2})}}%
\frac{dt}{t^{2}-1}}} \\
{} &=&{{}-I}\left( T,A\right) +{I}\left( T,B\right) {-I}\left( T,C\right) +{I%
}\left( T,0\right) \text{,}
\end{eqnarray*}
where
\begin{equation*}
{I}\left( T,A\right) ={\int\limits_{T}^{+\infty }{\log {\frac{t^{2}+A^{2}}{%
1+A^{2}}}\frac{dt}{t^{2}-1}}}
\end{equation*}
{Let $T=\tan \theta $. Then under substitution $t=\tan \tau $ we obtain}
\begin{eqnarray*}
{I}\left( T,A\right) {} &=& \\
&&\int\limits_{\theta }^{\frac{\pi }{2}}{\frac{{\log }\left( 1-\cos 2\tau
\cos 2\alpha \right) d\tau }{\cos 2\tau }}-\int\limits_{\theta }^{\frac{\pi
}{2}}{\frac{{\log }\left( 1+\cos 2\tau \right) d\tau }{\cos 2\tau }} \\
&=&\delta \left( \alpha ,\theta \right) -\delta \left( \frac{\pi }{2},\theta
\right) \text{.}
\end{eqnarray*}
Hence
\begin{equation*}
{\ }S(\alpha ,\beta ,\gamma )=-\delta (\alpha ,\theta )+\delta (\beta
,\theta )-\delta (\gamma ,\theta )+\delta (0,\theta )\text{.}
\end{equation*}
\end{proof}

\subsection{The relation $\protect\delta \left( \protect\alpha ,\protect%
\theta \right) $ to the {Lobachevskij }function}

Let
\begin{equation*}
\Delta \left( \alpha ,\theta \right) =-\int\limits_{\alpha -\theta }^{\alpha
+\theta }{\log }\left| {2\sin t}\right| dt=\Lambda \left( \alpha +\theta
\right) -\Lambda \left( \alpha -\theta \right)
\end{equation*}
and
\begin{equation*}
\delta \left( \alpha ,\theta \right) =\int\limits_{\theta }^{\frac{\pi }{2}}{%
\log }\left( {1-\cos 2\alpha \cos 2\tau }\right) \frac{d\tau }{\cos 2\tau }%
\text{.}
\end{equation*}

\begin{proposition}
Let $\alpha ,\,\theta ,\bigskip \widetilde{\theta }\in R$, $|
\widetilde{\theta }| $ $<\pi /4$ and $\tan \widetilde{\theta
}=\tanh \theta $. Then
\begin{equation*}
i\,\delta \left( \alpha ,\widetilde{\theta }\right) -i\,\delta \left( \alpha
,0\right) =\Delta \left( \alpha ,i\theta \right) -\Delta \left( \frac{\pi }{4%
},i\theta \right) \text{.}
\end{equation*}
\end{proposition}

\begin{proof}
Let
\begin{equation*}
F\left( \alpha ,\theta \right) =i\,\delta \left( \alpha ,\widetilde{\theta }%
\right) -i\,\delta \left( \alpha ,0\right) -\Delta \left( \alpha ,i\theta
\right) +\Delta \left( \frac{\pi }{4},i\theta \right) \text{.}
\end{equation*}
We have to show that $F\left( \alpha ,\theta \right) $ satisfy the equation
\begin{equation}
\frac{\partial F\left( \alpha ,\theta \right) }{\partial {\theta }}=0\text{,}
\label{deq}
\end{equation}
with the initial condition $F\left( \alpha ,0\right) =0$. It means $F\left(
\alpha ,\theta \right) =0$ and the proposition follows.

Note, that
\begin{equation}
\frac{\partial \widetilde{\theta }}{\partial \theta }=\frac{\cos ^{2}%
\widetilde{\theta }}{\cosh ^{2}\theta }=\frac{1-\tanh ^{2}\theta }{1+\tan
^{2}\widetilde{\theta }}=\frac{1-\tan ^{2}\widetilde{\theta }}{1+\tan ^{2}%
\widetilde{\theta }}=\cos 2\widetilde{\theta }\text{.}  \label{par}
\end{equation}
From $\left( \ref{par}\right) $ and definition of $\delta \left( \alpha ,{%
\theta }\right) $ we have
\begin{eqnarray}
\frac{\partial \left( i\,\delta \left( \alpha ,\widetilde{\theta }\right)
-i\,\delta \left( \alpha ,0\right) \right) }{\partial {\theta }} &=&i\frac{%
\partial \delta \left( \alpha ,\widetilde{\theta }\right) }{\partial
\widetilde{\theta }}\frac{\partial \widetilde{\theta }}{\partial {\theta }}
\label{der1} \\
&=&-i\log \left( 1-\cos 2\alpha \cos 2\widetilde{\theta }\right) \text{.}
\notag
\end{eqnarray}

From definition of $\Delta \left( \alpha ,\theta \right) $ we deduce
\begin{eqnarray*}
\frac{\partial \left( -\Delta \left( \alpha ,i\theta \right) +\Delta \left(
\frac{\pi }{4},i\theta \right) \right) }{\partial {\theta }} &=&i\log \left|
\frac{\sin \left( \alpha +i\theta \right) \sin \left( \alpha -i\theta
\right) }{\sin \left( \frac{\pi }{4}+i\theta \right) \sin \left( \frac{\pi }{%
4}-i\theta \right) }\right|  \\
&=&i\log \left| \frac{\cos 2i\theta -\cos 2\alpha }{\cos 2i\theta -\cos
\frac{\pi }{2}}\right|  \\
&=&i\log \left| 1-\frac{\cos 2\alpha }{\cos 2i\theta }\right| =i\log \left(
1-\frac{\cos 2\alpha }{\cosh 2\theta }\right) \text{.}
\end{eqnarray*}
Since
\begin{equation*}
\cos 2\widetilde{\theta }=\frac{1-\tan ^{2}\widetilde{\theta }}{1+\tan ^{2}%
\widetilde{\theta }}=\frac{1-\tanh ^{2}\theta }{1+\tanh ^{2}\theta }=\frac{1%
}{\cosh 2\theta }
\end{equation*}
we obtain
\begin{equation}
\frac{\partial \left( -\Delta \left( \alpha ,i\theta \right) +\Delta \left(
\frac{\pi }{4},i\theta \right) \right) }{\partial {\theta }}=i\log \left(
1-\cos 2\alpha \cos 2\widetilde{\theta }\right) \text{.}  \label{der2}
\end{equation}
The equalities $\left( \ref{der1}\right) $ and $\left( \ref{der2}\right) $
give $\left( \ref{deq}\right) $. The initial condition $F\left( \alpha
,0\right) =0$ follows directly from definitions of $\delta \left( \alpha
,\theta \right) $ and $\Delta \left( \alpha ,\theta \right) $.
\end{proof}

\subsection{\protect\bigskip The relation $\protect\delta \left( \protect%
\alpha ,\protect\theta \right) $ to the {Dilogarithm}}

Here we consider the relation
\begin{equation*}
\delta \left( \alpha ,\theta \right) =\int\limits_{\theta }^{\frac{\pi }{2}}{%
\log }\left( {1-\cos 2\alpha \cos 2\tau }\right) \frac{d\tau }{\cos 2\tau }
\end{equation*}
to the {Dilogarithm} function (\cite{Le}, p.~292)
\begin{equation*}
{\rm Li_{2} \,}\left( r,t\right) =-\frac{1}{2}\int\limits_{0}^{r}{\log }\left( {%
1-2x\cos t+x}^{2}\right) \frac{dx}{x}\text{.}
\end{equation*}

Recall that $\delta \left( \alpha ,\frac{\pi }{4}\right) $ is a $\pi $%
-periodic function with respect to $\alpha $ and
\begin{equation*}
\delta \left( \alpha ,\frac{\pi }{4}\right) =\left( \frac{\pi }{2}-\left|
\frac{\pi }{2}-\alpha \right| \right) ^{2}-\frac{\pi ^{2}}{16}
\end{equation*}
for $0<\alpha <\pi $.

\begin{proposition}
Let $\alpha \in R$ and $-\pi /4<\theta <3\pi /4$. Then
\begin{equation*}
\delta \left( \alpha ,\theta \right) -\delta \left( \alpha ,\frac{\pi }{4}%
\right) ={\rm Li_{2} \,}\left( \tan \left( \frac{\pi }{4}-\theta \right) ,\frac{\pi }{%
2}\right) -{\rm Li_{2} \,}\left( \tan \left( \frac{\pi }{4}-\theta
\right) ,2\alpha \right) \text{.}
\end{equation*}
\end{proposition}

\begin{proof}
Since
\begin{equation*}
{\rm Li_{2} \,}\left( \tan \mu ,\alpha \right) =-\int\limits_{0}^{\mu }{\log }\frac{%
1-\sin 2\mu \cos \alpha }{\cos ^{2}\mu }\frac{d\mu }{\sin 2\mu }
\end{equation*}
for $-\pi /2<\mu <\pi /2$ we have
\begin{gather}
{\rm Li_{2} \,}(\tan \mu ,\frac{\pi }{2})-{\rm Li_{2} \,}\left(
\tan \mu ,2\alpha \right) = \label{q1} \\
\int\limits_{0}^{\mu }{\log }\left( 1-\sin 2\mu \cos 2\alpha \right) \frac{%
d\mu }{\sin 2\mu }=-\delta \left( \alpha ,\mu -\frac{\pi }{4}\right) +\delta
\left( \alpha ,-\frac{\pi }{4}\right) .  \notag
\end{gather}

The property $\left( ii\right) $ of function $\delta \left( \alpha ,\mu
\right) $ gives
\begin{equation}
\delta \left( \alpha ,\frac{\pi }{4}-\mu \right) -\delta \left( \alpha ,%
\frac{\pi }{4}\right) =-\delta \left( \alpha ,\mu -\frac{\pi }{4}\right)
+\delta \left( \alpha ,-\frac{\pi }{4}\right)   \label{q2}
\end{equation}
and under substitution $\theta =\pi /4-\mu $ we obtain from $\left( \ref{q1}%
\right) $ and $\left( \ref{q2}\right) $
\begin{equation*}
\delta \left( \alpha ,\theta \right) -\delta \left( \alpha ,\frac{\pi }{4}%
\right) ={\rm Li_{2} \,}\left( \tan \left( \frac{\pi }{4}-\theta \right) ,\frac{\pi }{%
2}\right) -{\rm Li_{2} \,}\left( \tan \left( \frac{\pi }{4}-\theta
\right) ,2\alpha \right)
\end{equation*}
for $-\pi /4<\theta <3\pi /4$.
\end{proof}

\end{document}